%% file: TuringJumpsViaProvability.tex
\documentclass{article}

\input{Preamble}

\usepackage{bm}
\usepackage[affil-it]{authblk}

\newcommand{\trueBox}[1]{[#1]^{\sf True}}
\newcommand{\trueDiamond}[1]{\la #1 \ra^{\sf True}}

\newcommand{\boxBox}[1]{[#1]^{\Box}}
\newcommand{\boxDiamond}[1]{\la #1 \ra^{\Box}}

\newcommand{\omegaBox}[1]{[#1]^{\sf Omega}}
\newcommand{\omegaDiamond}[1]{\la #1 \ra^{\sf Omega}}

\begin{document}

\title{Turing jumps through provability}
\author{Joost J. Joosten\thanks{jjoosten@ub.edu}}
\affil{University of Barcelona}

\maketitle
\begin{abstract}
Fixing some computably enumerable theory $T$, the Friedman-Goldfarb-Harrington (FGH) theorem says that over elementary arithmetic, each $\Sigma_1$ formula is equivalent to some formula of the form $\Box_T \varphi$ provided that $T$ is consistent. In this paper we give various generalizations of the FGH theorem. In particular, for $n>1$ we relate $\Sigma_{n}$ formulas to provability statements $\trueBox{n}_T\varphi$ which are a formalization of ``provable in $T$ together with all true $\Sigma_{n+1}$ sentences". As a corollary we conclude that each $\trueBox{n}_T$ is $\Sigma_{n+1}$-complete. 

This observation yields us to consider a recursively defined hierarchy of provability predicates $[n+1]^\Box_T$ which look a lot like $\trueBox{n+1}_T$ except that where $\trueBox{n+1}_T$ calls upon the oracle of all true $\Sigma_{n+2}$ sentences, the $[n+1]^\Box_T$ recursively calls upon the oracle of all true sentences of the form $\la n \ra_T^\Box\phi$. As such we obtain a `syntax-light' characterization of $\Sigma_{n+1}$ definability whence of Turing jumps which is readily extended beyond the finite. Moreover, we observe that the corresponding provability predicates $[n+1]_T^\Box$ are well behaved in that together they provide a sound interpretation of the polymodal provability logic $\glp_\omega$.
\end{abstract}

\section{Introduction}

In first order arithmetic we have natural syntactical definitions that correspond to finite iterations of the Turing jump. Recall that a sentence in first order logic in the standard language of arithmetic is $\Sigma_{n+1}$ if it starts with a block of alternating quantifiers of length $n+1$ where the leftmost quantifier is existential and where the block of quantifiers is followed by a decidable formula only containing bounded quantification. There are various results known that relate these formula classes to computational complexity classes.

For example, a set of natural numbers is many-one reducible to the $n$-th Turing jump of the empty set if and only if it is $\Sigma_n$ definable (on the standard model of the natural numbers). Likewise, a set of natural numbers is Turing-reducible to the $n$-th Turing jump of the empty set if and only if it and its complement can be defined on the standard model of the natural numbers by a $\Sigma_{n+1}$ formula. Similarly, a set of natural numbers is computably enumerable relative to the $n$-th Turing jump of the empty set if and only if it can be defined by a $\Sigma_{n+1}$ formula.

In this paper we shall use the fact that various provability predicates are Turing complete in a certain sense so that we can give alternative characterizations for the finite Turing jumps. 
A central ingredient in proving our results come from generalizations of the so-called FGH Theorem.

The FGH Theorem (for Friedman-Goldfarb-Harrington) tells us that for any c.e.~theory $T$ we have provably in Elementary Arithmetic that each $\Sigma_1$ sentence $\sigma$ is equivalent to a provability statement of $T$, provided $T$ is consistent. In symbols, 
\[
\forall \, \sigma{\in}\Sigma_1 \, \exists \psi \ \  \ea \vdash \, \Diamond_T\top \to \big( \sigma \leftrightarrow \Box_T \psi \big).
\]
Here, as usual $\Box_T$ denotes a natural formalization of provability in $T$ and the $\Diamond_T$ stands for the dual consistency assertion. In this paper we give various generalizations of the FGH theorem. In particular we prove that the theorem holds for the provability notion $[n]_T$: provable in $T$ together with all true $\Sigma_{n+1}$ formulas. As a corollary we conclude that each $\trueBox{n}_T$ is $\Sigma_{n+1}$-complete. 

This observation yields us to consider a recursively defined hierarchy of provability predicates $[n+1]^\Box_T$ which look a lot like $\trueBox{n+1}_T$ except that where $\trueBox{n+1}_T$ calls upon the oracle of all true $\Sigma_{n+2}$ sentences, the $[n+1]^\Box_T$ recursively calls upon the oracle of all true sentences of the form $\la n \ra_T^\Box\phi$. 

As such we obtain a `syntax-light' characterization of $\Sigma_{n+1}$ definability whence of Turing jumps which is readily extended beyond the finite. Moreover, we observe that the corresponding provability predicates $[n+1]_T^\Box$ are well behaved in that together they provide a sound interpretation of the polymodal provability logic $\glp_\omega$.

\section{Preliminaries}

We shall work with theories with identity in the language $\{  0,1, \exp, +, \cdot, < \}$ of arithmetic where $\exp$ denotes the unary function $x\mapsto 2^x$. We define $\Delta_0=\Sigma_0=\Pi_0$ formulas as those where all quantifiers occur bounded, that is, we only allow quantifiers of the form $\forall\,  x{<}t$ or $\exists\, x{<}t$ where $t$ is some term not containing $x$. We inductively define $\Sigma_n, \Pi_n \subset \Pi_{n+1}$ and $\Sigma_n, \Pi_n \subset \Sigma_{n+1}$; if $\phi, \psi \in \Pi_{n+1}$, then $\forall x \ \phi, \phi \wedge \psi, \phi \vee \psi \in \Pi_{n+1}$ and likewise, if $\phi, \psi \in \Sigma_{n+1}$, then $\exists x \ \phi, \phi \wedge \psi, \phi \vee \psi \in \Sigma_{n+1}$.

We shall write $\Sigma_{n+1}!$ for formulas $\varphi$ of the form $\exists x\ \varphi_0$ with $\varphi_0 \in \Pi_n$. We will work in the absence of strong versions of (bounded) collection ${\sf B}_\varphi$ which is defined as 
\[
{\sf B}_\varphi \ := \ \forall z\, \forall \vec u \ \Big( \forall \, x{<}z \, \exists y \ \varphi(x,y, \vec u) \to \exists y' \,\forall \, x{<}z \, \exists \, y{<}y'\ \varphi(x,y, \vec u) \Big).
\]
Therefore, we will consider the formula class $\Sigma_{n+1,1}$ consisting of existentially quantified disjunctions and conjunctions of $\Sigma_{n+1}$ formulas with bounded quantifiers over them. To be more precise, we first inductively define 
\[
\Sigma_{n+1,b}:= \Sigma_{n+1}\mid(\Sigma_{n+1,b}\circ \Sigma_{n+1,b})\mid (\mathcal Q\, x{<}y \ \Sigma_{n+1,b})
\]
with $\circ \in \{ \wedge, \vee\}$ and $\mathcal Q \in \{ \forall, \exists\}$. Next we define the $\Sigma_{n+1,1}$ formulas to be of the form $\exists x \ \phi$ with $\phi \in \Sigma_{n+1,b}$.



The theory of \emph{elementary arithmetic}, \ea, is axiomatized by the defining axioms for $\{  0,1, \exp, +, \cdot, < \}$ together with induction for all $\Delta_0$ formulas. The theory \emph{Peano Arithmetic}, \pa, is as \ea but now allowing induction axioms for any first order formula. 

It is well known that \pa proves any instance ${\sf B}_\varphi$ of collection so that in particular each $\Sigma_{n+1,1}$ sentence is equivalent to some $\Sigma_{n+1}$ sentence.
Clearly we have that $\Sigma_{n+1}! \subset \Sigma_{n+1}$. Using coding techniques, it is clear that each $\Sigma_{n+1}$ formula is within \ea equivalent to a $\Sigma_{n+1}!$ formula. 

For us, a computably enumerable (c.e.) theory $T$ is understood to be given by a $\Delta_0$ formula that defines the set of codes of the primitive recursive set of axioms of $T$. We will employ standard formalizations of meta-mathematical properties like ${\tt Proof}_T(x,y)$ for ``$x$ is the G\"odel number of a proof from the axioms of $T$ of the formula whose G\"odel number is $y$". We shall often refrain from distinguishing a syntactical object $\varphi$ from its G\"odel number $\ulcorner \varphi \urcorner$ or from a syntactical representation of its G\"odel number. 

We will write $\Box_T \varphi$ for the $\Sigma_1!$ formula $\exists x\ {\tt Proof}_T(x,\varphi)$ and $\Diamond_T\varphi$ for $\neg \Box_T \neg \varphi$. By $\Box_T \varphi(\dot x)$ we will denote a formula which contains the free variable $x$, that expresses that for each value of $x$ the formula $\varphi(\overline x)$ is provable in $T$. Here, $\overline x$ denotes a syntactical representation of the number $x$.

By $\Sigma_1$ completeness we refer to the fact that for any true $\Sigma_1$ sentence $\sigma$ we have that $\ea \vdash \sigma$. It is well-known that \ea proves a formalized version of this: for any $\Sigma_1$ formula $\sigma(x)$ and any c.e.~theory $T$ we have $\ea \vdash \sigma(x) \to \Box_T\sigma(\dot x)$.

\section{The FGH theorem and generalizations}

In this section we shall be dealing with various so-called \emph{witness-comparison arguments} where the order of (least) witnesses to existential sentences is important. The first and most emblematic such argument occurred in the proof of Rosser's theorem which is a strengthening of G\"odel's first incompleteness theorem.

\begin{theorem}[Rosser's Theorem]\label{theorem:RossersTheorem}
Let $T$ be a consistent c.e.~theory extending \ea. There is some $\rho \in \Sigma_1$ which is undecidable in $T$. That is,
\[
\begin{array}{lr}
T \nvdash \rho & \mbox{\ \ \ and,} \\
T \nvdash \neg \rho .
\end{array}
\]
\end{theorem}

For rhetoric reasons we shall below include a standard proof of this celebrated result. Before doing so, we first need some notation. 


\begin{definition}
For $\phi := \exists x\  \phi_0 (x)$ and $\psi := \exists x\  \psi_0(x)$ we define 
\[
\begin{array}{llll}
\phi \leq \psi & := & \exists x \ \big(\phi_0(x) \wedge \forall \, y{<}x_0 \neg \psi_0(y)\big) & \mbox{ and, }\\
\phi < \psi & := & \exists x \ \big(\phi_0(x) \wedge \forall \, y{\leq}x \ \neg \psi_0(x)\big).& \\
\end{array}
\]
\end{definition}

Statements of the form $\phi \leq \psi$ or $\phi < \psi$ with $\phi, \psi \in \Sigma_{n+1}$ are called witness-comparison statements. Let us now collect some easy principles about witness-comparison statements whose elementary proofs we leave as an exercise.

\begin{lemma}\label{theorem:BasicPropertiesWitnesscomparisons}
For $A$ and $B$ in $\Sigma_{n+1}$ we have
\begin{enumerate}
\item
$\ea \vdash (A < B) \to (A \leq B)$;

\item
$\ea \vdash (A<B) \wedge (B\leq C) \to (A<C)$;

\item
$\ea \vdash (A\leq B) \wedge (B < C) \to (A<C)$;

\item
$\ea \vdash (A\leq B) \wedge (B\leq C) \to (A\leq C)$;

\item
$\ea \vdash (A \leq B) \to \neg (B< A)$ and consequently;

\item 
$\ea \vdash (A < B) \to \neg (B \leq A)$\label{item:AbelowBImpliesNotBatMostA:theorem:BasicPropertiesWitnesscomparisons};

\item
$\ea \vdash [(B\leq B) \vee (A\leq A)] \ \to \ [(A\leq B) \vee (B<A)]$; \label{item:sometimesLinearOrder:theorem:BasicPropertiesWitnesscomparisons}

\item
$\ea \vdash (A \leq B ) \to A$\label{item:PositiveInfo:theorem:BasicPropertiesWitnesscomparisons};

\item
$\ea \vdash A \wedge \neg B \to (A<B)$ \label{item:AandNotB:theorem:BasicPropertiesWitnesscomparisons};

\item \label{item:newitem:theorem:BasicPropertiesWitnesscomparisons}
$\ea \vdash A \wedge \neg (A\leq B) \to B$.

\item
Both $C< D$ and $C\leq D$ are of complexity $\Sigma_{n+1,1}$ if $C, D{\in} \Sigma_{n+1,1}$.

\end{enumerate}
\end{lemma}

We can now present a concise proof of Rosser's theorem.

\begin{proof}
We consider a fixpoint $\rho$ so that $T\vdash \rho \leftrightarrow (\Box_T \neg \rho \leq \Box_T \rho)$. 

If $T\vdash \rho$, then for some number $n$ we have ${\tt Proof}_T(n, \rho)$. Since $T$ is consistent we also have $\forall \, m{\leq}n \neg {\tt Proof}_T(m,\neg \rho)$. Thus, by $\Sigma_1$ completeness we have $T\vdash \Box_T \rho < \Box_T \neg \rho$ whence $T\vdash \neg \rho$; a contradiction.

Likewise, if $T\vdash \neg \rho$ we may conclude $\Box_T \neg \rho \leq \Box_T \rho$ so that $T\vdash \rho$.
\end{proof}

We would like to stress that it is actually quite remarkable that witness comparison arguments on statements involving G\"odel numbering can be used to yield any sensible information at all, since by tweaking the G\"odel numbering in a primitive recursive fashion we can always flip the order of the codes of any two syntactical objects. Of course, as we use fixpoints, after tweaking the G\"odel numbering, the corresponding fixpoint also changes. But still, it is remarkable that the useful witness comparison information of the fixpoint cannot be destroyed by tweaking the underlying G\"odel numbering.

We now turn our attention to another theorem, a proof of which can succinctly be given using witness comparison arguments: the \emph{FGH theorem}. The initials FGH refer to \emph{Friedman}, \emph{Goldfarb} and \emph{Harrington} who all substantially contributed to the theorem and we refer to \cite{Visser:2005:FaithAndFalsity} for historical details. 

Basically, the FGH theorem says that given any c.e.~theory $T$, any $\Sigma_1$ sentence is provably equivalent to a provability statement of the form $\Box_T \varphi$, modulo the consistency of $T$. The proof we give here is a slight modification of the one presented in \cite{Visser:2005:FaithAndFalsity}. The most important improvement is that we avoid the use of the least-number principle so that the proof becomes amenable for generalizations without a need to increase the strength of the base theory. 

\begin{theorem}[FGH theorem]\label{theorem:FGHTheorem}
Let $T$  be any computably enumerable theory extending \ea. For each $\sigma \in \Sigma_1$ we have that there is some $\rho \in \Sigma_1$ so that 
\[
\ea \vdash \Diamond_T\top \to \ \big (\sigma  \leftrightarrow \Box_T \rho \big ).
\]
\end{theorem} 

\begin{proof}
As in \cite{Visser:2005:FaithAndFalsity} we consider the fixpoint $\rho \in \Sigma_1$ for which $\ea \vdash \rho \,  \leftrightarrow \, (\sigma \leq \Box_T \rho)$. 
Without loss of generality we may assume that $\sigma\in \Sigma_1!$ so that both $\sigma \leq \Box_T \rho$ and $\Box_T \rho< \sigma$ are $\Sigma_1$. We now reason in \ea, assume $\Diamond_T \top$ and set out to prove $\sigma  \leftrightarrow \Box_T \rho$.

$(\to)$: assume for a contradiction that $\sigma$ and $\neg \Box_T \rho$. By Lemma \ref{theorem:BasicPropertiesWitnesscomparisons}.\ref{item:AandNotB:theorem:BasicPropertiesWitnesscomparisons} we conclude $\sigma \leq \Box_T \rho$, i.e., $\rho$. By provable $\Sigma_1$ completeness we get $\Box_T \rho$.

$(\leftarrow)$: assume for a contradiction that $\neg \sigma$ and $\Box_T \rho$. Again, we conclude $\Box_T \rho < \sigma$ so that $\Box_T (\Box_T \rho < \sigma)$ whence $\Box_T \neg \rho$ so that $\Box_T \bot$ contradicting the assumption $\Diamond_T \top$.
\end{proof}

A very special feature of $\rho$ from the above proof is that it is of complexity $\Sigma_1$ and that $\neg \rho$ is implied a by a related $\Sigma_1$ formula. Thus, we clearly provably have $\rho \to \Box_T \rho$ but in general we do not have $\neg \rho \to \Box_T \neg \rho$. But, due to the nature of $\rho$ we have that $\Box_T \neg \rho$ follows from the $\Sigma_1$ statement that is slightly stronger than $\neg \rho$, namely $\Box_T \rho < \sigma$: since provably $(\Box_T \rho < \sigma) \to \Box_T(\Box_T \rho < \sigma)$ and
\[
\begin{array}{lll}
(\Box_T \rho < \sigma) & \to & \neg (\sigma \leq \Box_T \rho)\\
 & \to & \neg \rho .\\
\end{array}
\]
The least number principle for $\Sigma_n$ formulas, ${\sf L}\Sigma_n$, says that for any $\psi \in \Sigma_n$ we have $\exists x\, \psi(x) \to \exists x\, (\psi(x) \wedge \forall\, y{<}x\ \neg \psi(y))$.  Of course, using ${\sf L}\Sigma_0$ and by Lemma \ref{theorem:BasicPropertiesWitnesscomparisons}.\ref{item:sometimesLinearOrder:theorem:BasicPropertiesWitnesscomparisons}, $\neg \rho$ and $\Box \rho < \sigma$ are provably equivalent under the assumption that $\Box \rho \vee \sigma$.

%
%

%

We will be interested in generalizing the FGH theorem to $\Sigma_n$ formulas using ever stronger notions of provability. Visser's proof of the FGH theorem as presented in \cite{Visser:2005:FaithAndFalsity} used an application of the least number principle for $\Delta_0$ formulas in the guise of $A \to (A\leq A)$. Thus, a direct generalization of Visser's argument to stronger provability notions would call for stronger and stronger arithmetical principles:

\begin{lemma}
The schema $A \to (A\leq A)$ for $A \in \Sigma_{n+1}!$ is over \ea provably equivalent to the least-number principle for $\Pi_n$ formulas. 
\end{lemma}

However, since our proof of the FGH theorem did not use the minimal number principle, we shall now see how the above argument generalizes to other provability predicates. By $\trueBox{n+1}_T$ we will denote the formalization of the predicate ``provable in $T$ together with all true $\Sigma_{n+2}$ sentences''. For convenience, we set $[0]_T := \Box_T$. Basically, for $n>0$, the predicate $\trueBox{n}_T \varphi$ will be a formalization of ``there is a sequence $\pi_0,\ldots, \pi_m$ so that each $\pi_i$ is either an axiom of $T$, or a true $\Sigma_{n+1}$ sentence, or a propositional logical tautology, or a consequence of some rule of $T$ using earlier elements in the sequence as antecedents". Thus, it is clear that for recursive theories $T$ we can write $\trueBox{n}_T$ by a $\Sigma_{n+1,1}$-formula. Also, we have provable $\Sigma_{n+1,1}$ completeness for these predicates, that is:

\begin{lemma}\label{theorem:SigmaNplusOneCompleteness}
Let $T$ be a c.e.~theory extending \ea and let $\phi$ be a $\Sigma_{n+1,1}$ formula. We have that
\[
\ea \vdash \phi(x) \to \trueBox{n}_T \phi(\dot x).
\]
\end{lemma}

\begin{proof}
Given $\phi (z,y_1,\ldots, y_k,x) \in \Sigma_{n+1}$, reason in \ea, fix $x_1,\ldots, x_k, x$ and, assume 
\[
\exists z\, \mathcal Q_1\, y_1{<}x_1 \ldots \mathcal Q_k \, y_k{<}x_k  \ \phi (z, y_1,\ldots, y_k,x)
\]
where $\mathcal Q_1\, y_1{<}x_1 \ldots \mathcal Q_k \, y_k{<}x_k $ is some block of bounded quantifiers. Thus, for some $a$ we have $\mathcal Q_1\, y_1{<}x_1 \ldots \mathcal Q_k \, y_k{<}x_k  \ \phi (a, y_1,\ldots, y_k, x)$. Under the box, we can now replace each $\forall\, y_i{<}\overline {x_i}$ by $\bigwedge_{\overline{y_i}<\overline{x_i}}$ and each $\exists\, y_i{<}\overline{x_i}$ by $\bigvee_{\overline{y_i}<\overline{x_i}}$ so that by applying distributivity we see that $\mathcal Q_1\, y_1{<}\overline{x_1} \ldots \mathcal Q_k \, y_k{<}\overline{x_k}  \ \phi (\overline{a}, y_1,\ldots, y_k, \overline{x})$ is equivalent to a disjunctive normal form of bounded substitution instances of $\phi (\overline{a}, y_1,\ldots, y_k, \overline{x})$. Note that this operation in available within \ea since it only requires the totality of exponentiation. 

Outside the box we know that for some of these big conjunctions of bounded substitution instances of $\phi (a, y_1,\ldots, y_k,x)$ actually all of the conjuncts are true. Since each of those conjuncts is a true $\Sigma_{n+1}$ sentence, each conjunct is an axiom whence holds under the box. Thus, the whole conjunct is provable under the box whereby we obtain $\trueBox{n}_T \mathcal Q_1\, y_1{<}\overline{x_1} \ldots \mathcal Q_k \, y_k{<}\overline{x_k}  \ \phi (\overline{a}, y_1,\ldots, y_k, \overline{x})$ whence $\trueBox{n}_T \exists z\, \mathcal Q_1\, y_1{<}\overline{x_1} \ldots \mathcal Q_k \, y_k{<}\overline{x_k}  \ \phi (z, y_1,\ldots, y_k, \overline{x})$ as was to be shown. It is clear that this case suffices for the more general form of $\Sigma_{n+1,1}$ formulas.
\end{proof}

It is easy to check that the predicate $\trueBox{n}_T$ is well behaved. In particular one can check that all the axioms of the standard provability logic \gl as defined in the last section hold for it. Over \ea, the notion of $\trueBox{n}_T$ can be related to regular provability $\Box_T$ by the lemma below. 

By ${\tt FinSeq}(f)$ we denote a predicate that only holds on numbers that are codes of a finite sequence of G\"odel numbers and by $|f|$ we denote the length of such a sequence. Moreover, $f_i$ will denote the $i$th element of such a sequence $f$. With ${\tt True}_{\Sigma_{n+1}}$ we will denote a partial truth predicate for $\Sigma_{n+1}$ formulas and ${\tt Tr}_{\Sigma_{n+1}}$ will denote the set of true $\Sigma_{n+1}$ sentences.

\begin{lemma}
For any c.e.~theory $T$, we have that 
\[
\ea \vdash \trueBox{n}_T \varphi \leftrightarrow \exists f\, \Big( {\tt FinSeq}(f) \wedge \forall \, i {<}|f| \ {\tt True}_{\Sigma_{n+1}}(f_i) \, \wedge\, \Box_T \big( (\wedge_{i<|f|} f_i) \to \varphi\big)  \Big).
\]
\end{lemma}

\begin{proof}
We reason in \ea fixing some $\varphi$. The $\leftarrow$ direction follows directly from the formalized deduction theorem. For the other direction fix some $p$ with  ${\tt Proof}_{T+{\tt Tr}_{\Sigma_{n+1}}}(p, \varphi)$. We can express that $p_i$ is a true $\Sigma_{n+1}$ sentence in a $\Delta_0$ fashion simply by saying that it is not a propositional tautology, nor an axiom of $T$, nor the result of applying a rule to earlier elements in the sequence. Thus, by $\Delta_0$ induction on the length of $p$ we can prove that there is a sequence that collects all the true $\Sigma_{n+1}$ sentences.
\end{proof}

Note that our definition of $\trueBox{n}_T$ is slightly non-standard since in the literature (e.g.~\cite{Beklemishev:2004:ProvabilityAlgebrasAndOrdinals}) it is more common to define $\trueBox{n}_T$ using a $\Pi_n$ oracle rather than a $\Sigma_{n+1}$ oracle. With a $\Pi_n$ oracle one gets provable $\Sigma_{n+1}$ completeness but in the absence of ${\sf B}\Sigma_{n+1}$ not necessarily provable $\Sigma_{n+1,1}$ completeness. With our definition of $\trueBox{n}_T$, since $A\leq B \in \Sigma_{n+1,1}$ for $A,B \in \Sigma_{n+1,1}$, and since we provided a proof where the minimal number principle is avoided, the FGH theorem smoothly generalizes to the new setting. 

\begin{theorem}\label{theorem:FGHforTrueBox}
Let $T$  be any computably enumerable theory extending \ea and let $n<\omega$. For each $\sigma \in \Sigma_{n+1,1}$ we have that there is some $\rho_n \in \Sigma_{n+1,1}$ so that 
\[
\ea \vdash \trueDiamond{n}_T \top \to \ \big(\sigma  \leftrightarrow \trueBox{n}_T \rho_n\big).
\]
\end{theorem}

\begin{proof}
The proof runs entirely analogue to the proof of Theorem \ref{theorem:FGHTheorem}. Thus, for each number $n$ we consider the fixpoint $\rho_n$ so that $\ea \vdash \rho_n \leftrightarrow (\sigma \leq \trueBox{n}_T \rho_n)$. Note that both $\rho_n$ and $\trueBox{n}_T \rho < \sigma$ are $\Sigma_{n+1,1}$ whence by Lemma \ref{theorem:SigmaNplusOneCompleteness} we can apply provable $\Sigma_{n+1,1}$ completeness to them.
\end{proof}
As an easy corollary we get that $\trueBox{n}_T$ formulas are closed not only under conjunction, as is well know, but also under disjunctions.
Note that the FGH theorem yields that provably $\trueDiamond{n}_T \top \to (\sigma \leftrightarrow \trueBox{n}_T\rho_n)$. Using the propositional tautology
\[
( \neg A \to C) \to \Big[ \big ( A \to (B\leftrightarrow C) \big) \ \leftrightarrow \ \big ( (\neg A \vee B) \leftrightarrow C \big) \Big]
\]
and $\trueBox{n}_T \bot \to \trueBox{n}_T \rho_n$ we see that this is equivalent to $\big (\trueBox{n}_T \bot \vee \sigma \big) \leftrightarrow \trueBox{n}_T \rho_n$.
\begin{corollary}\label{theorem:finiteTrueBoxFormulasClosedUnderDisjunctions}
Let $T$ be a c.e.~theory extending \ea and let $n\in \mathbb N$. For each formulas $\varphi, \psi$ there is some $\sigma \in \Sigma_{n+1}$ so that 
\[
T\vdash (\trueBox{n}_T \varphi \ \vee \trueBox{n}_T \psi) \ \leftrightarrow \ \trueBox{n}_T \sigma .
\]
\end{corollary}
\begin{proof}
We consider some $\sigma \in \Sigma_{n+1,1}$ so that provably $\sigma \leftrightarrow (\trueBox{n}_T \varphi \ \vee \trueBox{n}_T \psi)$.
By Theorem \ref{theorem:FGHforTrueBox} applied to this $\sigma$ we provably have that $(\trueBox{n}_T \varphi \ \vee \trueBox{n}_T \psi)\vee \trueBox{n}_T \bot$ is equivalent to $\trueBox{n}_T \varphi \ \vee \trueBox{n}_T \psi$.
\end{proof}

As another corollary of Theorem \ref{theorem:FGHforTrueBox} we see that, in a sense, the notion of $n$-provability is $\Sigma_{n+1}$ complete. To establish this, we will first need a particular version of the fixpoint lemma.

\begin{lemma}\label{theorem:simpleDynamicFixpointLemma}
Let $\psi (x,y)$ be a formula whose free variables are amongst $\{ x,y\}$. There is a formula $\varphi (y)$ so that 
\[
\ea \vdash \varphi (y) \leftrightarrow \psi (\ulcorner \varphi(\dot y)\urcorner, y).
\]
\end{lemma}

\begin{proof}
The lemma easily follows from the well-known version of the fixpoint lemma by which for each formula $\psi (x,y)$ there is a formula $\varphi (y)$ so that 
\begin{equation}\label{equation:fixpointLemmaSimple}
\ea \vdash \varphi (y) \leftrightarrow \psi (\ulcorner \varphi(y)\urcorner, y)
\end{equation}
(e.g., the \emph{generalized diagonal lemma} from Boolos' \cite{Boolos:1993:LogicOfProvability}, Chapter 3).

By $z = {\sf sub}(u,v,w)$ we denote the formula that expresses that $z$ is the G\"odel number of the result of substituting the numeral of $w$ for the variable whose G\"odel number is $v$ into the formula whose G\"odel number is $u$. By $\psi(x/z,y)$ we denote the result of substituting $z$ for $x$ in $\psi(x,y)$. We now consider the formula 
\[
\exists z\ \big(z= {\sf sub}(x, \ulcorner y \urcorner, y) \ \wedge \ \psi(x/z,y)\big)
\]
and apply Equation \eqref{equation:fixpointLemmaSimple} to it to obtain the required fixpoint. Note that 
\[
\ea \vdash \Big[ \exists z\ \big(z= {\sf sub}(\ulcorner \varphi (y) \urcorner, \ulcorner y \urcorner, y) \ \wedge \ \psi(x/z,y)\big)\Big] \ \ \longleftrightarrow \ \  \psi(\ulcorner \varphi (\dot y) \urcorner,y).
\]
\end{proof}

With this lemma we can now prove $\Sigma_{n+1}$-completeness of the $\trueBox{n}_T$ provability predicate.


\begin{lemma}\label{theorem:TrueProvabilityTiesUpWithArithmeticalHierarchy}
Let $T$ be any sound c.e.~theory and let $A\subseteq \mathbb N$. The following are equivalent
\begin{enumerate}

\item\label{item:TuringReducible:theorem:TrueProvabilityTiesUpWithArithmeticalHierarchy}
$A$ is c.e.~in $\emptyset^{(n)}$;

\item\label{item:OneToOneReducible:theorem:TrueProvabilityTiesUpWithArithmeticalHierarchy}
$A$ is many-one reducible to $\emptyset^{(n+1)}$;

\item\label{item:Definable:theorem:TrueProvabilityTiesUpWithArithmeticalHierarchy}
$A$ is definable on the standard model by a $\Sigma_{n+1}$ formula;

\item\label{item:BoxDefinable:theorem:TrueProvabilityTiesUpWithArithmeticalHierarchy}
$A$ is definable on the standard model by a formula of the form $\trueBox{n}_T \rho(\dot x)$;

\item\label{item:ReducedBoxDefinable:theorem:TrueProvabilityTiesUpWithArithmeticalHierarchy}
$A$ is definable on the standard model by a formula of the form $\trueBox{n}_T \rho(\dot x)$ where $\rho(x) \in \Sigma_{n+1,1}$;
\end{enumerate}
\end{lemma}

\begin{proof}
The equivalence of 
\ref{item:TuringReducible:theorem:TrueProvabilityTiesUpWithArithmeticalHierarchy}, 
\ref{item:OneToOneReducible:theorem:TrueProvabilityTiesUpWithArithmeticalHierarchy}, and 
\ref{item:Definable:theorem:TrueProvabilityTiesUpWithArithmeticalHierarchy} is just Post's theorem. 

The implication \ref{item:ReducedBoxDefinable:theorem:TrueProvabilityTiesUpWithArithmeticalHierarchy} $\Rightarrow$ \ref{item:BoxDefinable:theorem:TrueProvabilityTiesUpWithArithmeticalHierarchy} is trivial, and implication \ref{item:BoxDefinable:theorem:TrueProvabilityTiesUpWithArithmeticalHierarchy} $\Rightarrow$ \ref{item:Definable:theorem:TrueProvabilityTiesUpWithArithmeticalHierarchy} 
holds in virtue of $\trueBox{n}_T$ being a $\Sigma_{n+1,1}$ predicate which on $\mathbb N$ is equivalent to some $\Sigma_{n+1}$ formula, so it suffices to prove \ref{item:Definable:theorem:TrueProvabilityTiesUpWithArithmeticalHierarchy} 
$\Rightarrow$ 
\ref{item:ReducedBoxDefinable:theorem:TrueProvabilityTiesUpWithArithmeticalHierarchy}.

Thus, let the number $n$ be fixed and, let $A$ be a set of natural numbers so that for some $\sigma(x) \in \Sigma_{n+1}$ we have $m\in A \ \Longleftrightarrow \ \mathbb N \models \sigma (m)$. Using lemma \ref{theorem:simpleDynamicFixpointLemma} we find $\rho_n(x) \in \Sigma_{n+1,1}$ so that 
\[
\ea \vdash \forall x\ \big(\ \rho_n(x) \ \leftrightarrow \ [\sigma (x) \leq \trueBox{n}_T \rho_n (\dot x)] \ \big).
\]
Reasoning in \ea, we pick an arbitrary $x$ and repeat the reasoning as in the proof of Theorem \ref{theorem:FGHforTrueBox} to see that 
\[
\ea \vdash \trueDiamond{n}_T \top \to  \forall x\ \Big( \ \sigma(x) \ \ \leftrightarrow \ \ \big(\trueBox{n}_T \rho_n(\dot x)\big)\ \Big).
\]
Since \ea is sound and by assumption of $T$ also being sound we have for each $n$ that $\mathbb N \models \trueDiamond{n}_T \top$, we may conclude that for any number $m$,
\[
\mathbb N \models \sigma (m) \ \ \  \Longleftrightarrow \ \ \ \mathbb N \models\  \trueBox{n}_T \rho_n(\overline m)
\]
which was to be proven.
\end{proof}


Feferman showed in \cite{Feferman:1957:DegreesOfUnsolvability} that for each unsolvable Turing degree $\bm d$ there is a theory $U$ so that the Turing degree of $\{ \ulcorner \varphi \urcorner \mid U\vdash \varphi\}$ is $\bm d$. However, the theories that Feferman considered were formulated in the language of identity and in particular did not contain arithmetic. It is not hard to see that for theories that do contain arithmetic we can only attain degrees that arise as Turing jumps. 

\begin{lemma}
Let $A\subseteq \mathbb N$ be definable on $\mathbb N$ by $\alpha(x)$ with Turing degree $\bm a$. Then, the
\[
\mbox{  Turing degree of  }\{  \psi \mid \ea + \{ \alpha(\overline n) \mid \mathbb N \models \alpha (n )\} \vdash \psi \} \mbox{ equals } {\bm a}'.
\]
%
%
\end{lemma}

\begin{proof}
Let us denote the Turing degree of a set $B$ by $\bm{|}B\bm{|}$. We see that 
\[
\bm{|}\{  \psi \mid \ea + \{ \alpha(\overline n) \mid \mathbb N \models \alpha (n )\} \vdash \psi \}\bm{|} \leq_T {\bm a}'
\]
since a Turing machine can enumerate all oracle proofs thereby reducing provability to the halting problem. Similarly, we see that 
\[
{\bm a}'  \leq_T \bm{|}\{  \psi \mid \ea + \{ \alpha(\overline n) \mid \mathbb N \models \alpha (n )\} \vdash \psi \}\bm{|}
\]
since we can code Turing machine computations in arithmetic.
\end{proof}

We thus see that the Turing degrees of theories that are defined by a minimal amount of arithmetic (\ea) together with some oracle, are entirely determined by the Turing degree of the corresponding oracle by means of the jump operator. By Friedberg's jump inversion theorem (\cite{Friedberg:1957:ACriterionForCompleteness}) we may thus conclude that any Turing degree above ${\bm 0}'$ can be attained as the decision problem of a theory containing arithmetic.  

In this light, Lemma \ref{theorem:TrueProvabilityTiesUpWithArithmeticalHierarchy} should not come as a surprise. However, in the lemma we have provided a natural subsequence of theories so that moreover, all the necessary reasoning for the reductions can be formalized in \ea.

Lemma \ref{theorem:TrueProvabilityTiesUpWithArithmeticalHierarchy} is stated entirely in terms of definability and computability but the proof tells us actually a bit more, namely that the FGH theorem is easily formalizable within \ea.

\begin{lemma}
For any c.e.~theory $T$ we have that 
\[
\ea \vdash \forall \, \sigma{\in}\Sigma_{n+1}\, \exists \, \rho{\in} \Sigma_{n+1} \ \Big ( \trueDiamond n_T \top \ \to \ \ \big(  {\sf True}_{\Sigma_{n+1}} (\sigma) \leftrightarrow \trueBox{n+1}_T \rho \big )\, \Big).
\]
\end{lemma}

\begin{proof}
One can simply formalize the proof of Theorem \ref{theorem:FGHforTrueBox} which is easy since mapping the G\"odel number of a $\Sigma_{n+1}$ sentence $\sigma$ to the G\"odel number of its  corresponding fixpoint $\rho_{n}$ is elementary. 

Alternatively, using Lemma \ref{theorem:simpleDynamicFixpointLemma} we find $\rho_n(x) \in \Sigma_{n+1}$ so that 
\[
\ea \vdash \forall x\ \big(\ \rho_n(x) \ \leftrightarrow \ [{\sf True}_{\Sigma_{n+1}}(x) \leq \trueBox{n}_T \rho_n (\dot x)] \ \big).
\]
so that by reasoning as in the proof of 
Lemma \ref{theorem:TrueProvabilityTiesUpWithArithmeticalHierarchy} we see that 
\[
\ea \vdash \trueDiamond{n}_T \top \ \to\ \forall x\ \Big( \ {\sf True}_{\Sigma_{n+1}}(x) \ \leftrightarrow \ \ \big(\trueBox{n}_T \rho_n(\dot x)\big)\ \Big).
\]
\end{proof}

For other notions of provability we get similar generalizations of the FGH theorem. In particular, let $\omegaBox{n}_T$  denote the formalization of the predicate ``provable in $T$ using at most $n$ nestings  of the omega rule''. Following the recursive scheme $\omegaBox{0}_T \varphi := \Box_T \varphi$ and, 
\[
\omegaBox{n+1}_T \varphi := \exists \psi \ \Big ( \forall x\ \omegaBox{n}_T \psi (\dot x) \ \wedge \ \Box_T (\forall x\ \psi (x) \to \varphi)\Big)
\]
we see that for c.e.~theories $T$ we can write $\omegaBox{n}_T$ by a $\Sigma_{2n+1}$-formula. Also, we have provable $\Sigma_{2n+1}$ completeness for these predicates, that is:

\begin{proposition}\label{theorem:SigmaTwoNplusOneCompleteness}
Let $T$ be a computable theory extending \ea and let $\phi$ be a $\Sigma_{2n+1}$ formula. We have that
\[
\ea \vdash \phi \to \omegaBox{n}_T \phi.
\]
\end{proposition}

\begin{proof}
By an external induction on $n$ where each inductive step requires the application of an additional omega-rule.
\end{proof}

This proposition is the omega-rule analogue of provable $\Sigma_{n+1,1}$ completeness for the $\trueBox{n}_T$ predicate. As a corollary we get an FGH Theorem for omega-provability.

\begin{corollary}\label{theorem:FGHforOmegaBox}
Let $T$  be any sound computably enumerable theory extending \ea and let $n<\omega$. For each $\sigma \in \Sigma_{2n+1}$ we have that there is some $\rho_n \in \Sigma_{2n+1,1}$ so that 
\[
\pa \vdash \omegaDiamond{n}_T \top \to \ \big(\sigma  \leftrightarrow \omegaBox{n}_T \rho_n\big).
\]
\end{corollary}

We have formulated this corollary over \pa so that $\Sigma_{2n+1,1}$ sentences are provably equivalent to $\Sigma_{2n+1}$ sentences using collection. Consequently, we can now also prove a definability result for the $\omegaBox{n}_T$ predicate.

\begin{lemma}\label{theorem:OmegaProvabilityTiesUpWithArithmeticalHierarchy}
Let $T$ be any c.e.~theory, let $n$ be a natural number, and let $A\subseteq \mathbb N$. The following are equivalent
\begin{enumerate}
\item
$A$ is c.e.~in $\emptyset^{(2n)}$;

\item
$A$ is definable on the standard model by a $\Sigma_{2n+1}$ formula;

\item
$A$ is definable on the standard model by a formula of the form $\omegaBox{n}_T \rho(\dot x)$;

\item
$A$ is definable on the standard model by a formula of the form $\omegaBox{n}_T \rho(\dot x)$ where $\rho(x) \in \Sigma_{2n+1,1}$;
\end{enumerate}
\end{lemma}

\begin{proof}
The proof of this lemma is analogous to the proof of Lemma \ref{theorem:TrueProvabilityTiesUpWithArithmeticalHierarchy} if one substitutes $\Sigma_{n+1}$ by $\Sigma_{2n+1}$, \ea by \pa,  and $\trueBox{n}_T$ by $\omegaBox{n}_T$. 
\end{proof}

By comparing lemmas \ref{theorem:TrueProvabilityTiesUpWithArithmeticalHierarchy} and \ref{theorem:OmegaProvabilityTiesUpWithArithmeticalHierarchy} we see that in a sense the hierarchy of formulas of the form $\trueBox{n}_T\varphi$ is more fine-grained than the hierarchy of formulas of the form $\omegaBox{n}_T\varphi$. 

The positive feature of the latter hierarchy is that it is defined solely in terms of provability whereas the former needs to call upon partial truth predicates. As such the $\omegaBox{n}_T$ hierarchy is more amenable to Turing jumps beyond the finite where no clear-cut syntactical characterizations along the lines of Post's correspondence theorem are available (see \cite{FernandezJoosten:2013:OmegaRuleInterpretationGLP}). 

The down-side to the $\omegaBox{n}_T$ hierarchy is that it runs outline with the Turing-jump hierarchy. In a forth-coming paper we propose a transfinite progression of provability notions in second order arithmetic that takes the best of both worlds: it is defined purely in terms of provability as in \eqref{equation:boxBoxDefinition}, synchronizes with the Turing-jump hierarchy as in Theorem \ref{theorem:recursionIsFineForFiniteNumbers}, and can be transfinitely extended along any ordinal $\Xi$ definable in second order logic yielding for a large class of second order theories a sound interpretation of a well-behaved logic called $\glp_\Xi$ as in Theorem \ref{theorem:glpOmegaSoundUnderBoxBoxInterpretation}. \\
\medskip

We end this section with some remarks on the fixpoint used in the proof of the FGH theorem. By no means, this fixpoint is the only one that works. The minor change where we consider $\rho \leftrightarrow (\sigma < \Box \rho)$ works with almost the same proof. But we also have a `dual' version of our fixpoint which gives the desired result.

\begin{lemma}
Let $T$ be a c.e.~theory and let $\sigma \in \Sigma_1$. 
\[
\mbox{If } \ \ea \vdash \rho \ \leftrightarrow \ (\Box_T \neg \rho \leq \sigma) \ \ \mbox{ then } \ \ \ea \vdash \sigma \vee \Box_T \bot \ \leftrightarrow \ \Box_T \neg \rho.
\]
\end{lemma}

\begin{proof}
We reason in \ea assuming $\rho \ \leftrightarrow \ (\Box_T \neg \rho \leq \sigma)$. 

 $\rightarrow:$
 We suppose $\sigma$. If $\rho$, then $\Box_T \neg \rho \leq \sigma$ whence $\Box_T \neg \rho$. If $\neg \rho$, then 
\begin{equation}\label{equation:dualFixPointFGHEquation}
\neg (\Box_T \neg \rho \leq \sigma).
\end{equation}
But from our assumption $\sigma$ we get by the minimal number principle that $\sigma \leq \sigma$ so that $(\Box_T \neg \rho \leq \sigma) \vee (\sigma < \Box_T \neg \rho)$ whence by \eqref{equation:dualFixPointFGHEquation} we get $\sigma < \Box_T \neg \rho$. By provable $\Sigma_{1}$-completeness we get $\Box_T (\sigma < \Box_T \neg \rho)$ whence $\Box_T \neg \rho$.

 $\leftarrow:$ If $\Box_T \neg \rho$ and $\neg \sigma$, then $\Box_T \neg \rho \leq \sigma$ whence $\rho$ and by $\Sigma_1$ completeness $\Box_T \rho$ so that $\Box_T \bot$.
\end{proof}

Since this proof uses the minimal number principle, it is not amenable for generalizations to stronger provability notions over a weak base theory. The following fixpoint does allow such generalizations. Note that the two lemma's are very similar yet not necessarily equivalent since fixpoints involving the witness-comparison relation are in general not closed under substituting logical equivalents.

\begin{lemma}
Let $T$ be a c.e.~theory and let $\sigma \in \Sigma_1$. 
\[
\mbox{If } \ \ea \vdash \rho \ \leftrightarrow \  \neg (\Box_T \rho < \sigma) \ \ \mbox{ then } \ \ \ea \vdash \sigma \vee \Box_T \bot \ \leftrightarrow \ \Box_T \rho.
\]
\end{lemma}

\begin{proof}
Again, we reason in \ea now avoiding the minimal number principle.

$\to$: Suppose $\sigma$. If $\sigma \leq \Box_T \rho$, then $\Box_T (\sigma \leq \Box_T \rho)$ whence $\Box_T \neg (\Box_T \rho < \sigma)$ so that $\Box \rho$. If $\neg (\sigma \leq \Box_T \rho)$, then since $\sigma$ we have $\Box_T \rho$.

$\leftarrow$: Suppose $\Box_T \rho$ and $\neg \sigma$. Then, $\Box_T \rho < \sigma$ whence also $\Box_T (\Box_T \rho < \sigma)$. Thus $\Box_T \neg \rho$, which together with the assumption that $\Box_T \rho$ yields $\Box_T \bot$.
\end{proof}

\section{An algebraic perspective}

In this section we shall recast the generalized FGH theorem in algebraic terms. Let us fix some c.e.~theory $T$. The generalized FGH theorem tells us that $\Sigma_{n+1}$ sentences are almost equivalent to sentences of the form $\trueBox{n}_T \phi$: that is, $\forall \sigma \in \Sigma_{n+1}$ there is $\rho$ such that $\ea \vdash \sigma \vee \trueBox{n}_T \bot \ \leftrightarrow \ \trueBox{n}_T \rho$. Thus, `almost equivalent' refers to that we always have to take the disjunct $\trueBox{n}_T \bot$ along. This is due to the fact that $\trueBox{n}_T \bot$ is minimal among all sentences of the form $\trueBox{n}_T\phi$ hence we cannot find $\trueBox{n}_T\phi$-equivalents of $\Sigma_{n+1,1}$ sentences below $\trueBox{n}_T\bot$. 

It is well know that there are plenty of such sentences. For example, the structure of $\Sigma_1$ sentences within the Lindenbaum algebra of any consistent c.e.~theory is known to be densely ordered by the following folklore fact which can be readily proven by an application of Rosser's Theorem.

\begin{fact}
Let $T$ be a c.e.~theory. For any two $\Sigma_1$ sentences $\sigma_0, \sigma_2$ so that $T\vdash \sigma_0 \to \sigma_2$ but $T\nvdash \sigma_2 \to \sigma_0$, there is some $\sigma_1$ which lies strictly in between $\sigma_0$ and $\sigma_2$. That is, $T\vdash (\sigma_0 \to \sigma_1) \wedge (\sigma_1 \to \sigma_2)$ but $T\nvdash \sigma_2 \to \sigma_1$ and $T\nvdash \sigma_1 \to \sigma_0$.
\end{fact}

\begin{proof}
Since $T\nvdash \sigma_2 \to \sigma_0$, we know that $T + \neg \sigma_0 + \sigma_2$ is consistent and will consider the corresponding Rosser sentence $\rho \in \Sigma_{1}$ as in Theorem \ref{theorem:RossersTheorem}. Now, $\sigma_1:= \sigma_0 \vee (\rho \wedge \sigma_2)$ will do the job: Clearly $T \vdash (\sigma_0 \to \sigma_1) \wedge (\sigma_1 \to \sigma_2)$; also $T + \neg \sigma_0 + \sigma_2 \nvdash \rho$ implies $T\nvdash \sigma_2 \to \sigma_1$, and $T + \neg \sigma_0 + \sigma_2 \nvdash \neg \rho$ implies $T\nvdash \sigma_1 \to \sigma_0$.
\end{proof}

Before we recast the FGH theory in terms of algebras we need some additional notation. Let $\Sigma_{n+1,1}/{\equiv_T}$ denote the set of $T$-equivalence classes of $\Sigma_{n+1,1}$ sentences and likewise, let $[n]_T \mathcal F /{\equiv_T}$ denote the set of $T$-equivalence classes of sentences of the form $\trueBox{n}_T \phi$. Moreover, let $(\Sigma_{n+1}\upharpoonright [n]_T \bot)/{\equiv_T}$ denote the set of $T$-equivalence classes of $\Sigma_{n+1}$ sentences $\varphi$ for which $T\vdash \trueBox{n}_T \bot \to \varphi$. 

By Corollary \ref{theorem:finiteTrueBoxFormulasClosedUnderDisjunctions} we can define $\vee^{[n]}: ([n]_T \mathcal F \times [n]_T \mathcal F) \to [n]_T \mathcal F$ so that 
\[
\vee^{[n]}(\trueBox{n}_T \varphi, \trueBox{n}_T \psi) := \trueBox{n}_T \chi \ \mbox{ where } \  T\vdash (\trueBox{n}_T \varphi \vee \trueBox{n}_T \psi ) \leftrightarrow  \trueBox{n}_T \chi.
\]
This map $\vee^{[n]}$ naturally extends to $[n] \mathcal F/{\equiv_T}$ and substructures thereof. Likewise, we define $\wedge^{[n]} (\trueBox{n}_T \varphi, \trueBox{n}_T \psi):= \trueBox{n}_T (\varphi\wedge\psi)$. With this notation the FGH theorem can be reinterpreted in terms of algebras.

\begin{lemma}
For any c.e.~theory $T$, the algebras
\begin{enumerate}
\item
$\la (\Sigma_{n+1}\upharpoonright [n]_T \bot)/{\equiv_T}, \wedge, \vee, [n]_T \ra$;

\item
$\la [n]_T \mathcal F/{\equiv_T}, \wedge^{[n]}, \vee^{[n]}, [n]_T \ra$;

\end{enumerate}
define the same subalgebra of the Lindenbaum algebra of $T$.
\end{lemma}

\begin{proof}
This follows directly from the above considerations. We can map $\sigma$ to $\trueBox{n+1}_T \rho_{n+1}$ as in the FGH theorem to establish a natural map on different representations of the same equivalence class within the Lindenbaum algebra of $T$.
\end{proof}

\section{Graded provability via Turing jumps}

We shall now see how the FGH theorem can be used to define graded provability notions $\boxBox{n}_T$ for $n\in \mathbb N$ which are defined using only provability notions yet which are $\Sigma_{n+1}$ complete very much in the same way as Lemma \ref{theorem:TrueProvabilityTiesUpWithArithmeticalHierarchy} told us that the $\trueBox{n}_T$ predicates are $\Sigma_{n+1}$ complete.
Recall that on the natural numbers we have 
\[
\trueBox{n+1}_T \varphi \ \ \Leftrightarrow \ \ \exists \, \pi \ \Big({\sf True}_{\Pi_{n+1}}(\pi) \wedge \Box_T ({\sf True}_{\Pi_{n+1}}(\pi) \to \phi)\Big)
\]
where $\trueBox{0}_T$ is nothing but $\Box_T$. It is easy to see that this equivalence is provable within \pa.
The idea now is to replace true $\Pi_{n+1}$ sentences by consistency statements, that is, by sentences of the form $\langle n \rangle_T \varphi$. This replacement will be done in a recursive fashion. Thus we can consider the following recursive scheme.
\begin{equation}\label{equation:OldboxBoxDefinition}
\boxBox{\overline{0}}_T\phi := \Box_T \phi, \ \ \mbox{ and } \ \ \boxBox{\overline{n+1}}_T \phi \ := \ \Box_T \phi \ \vee \ \exists \, \psi \ \big(\boxDiamond{\overline{n}}_T \psi \ \wedge \ \Box (\boxDiamond{\overline{n}}_T \psi \to \phi)\big).
\end{equation}
For this recursive scheme, we can easily prove various desirable properties. However, for the sake of generalizing the definition to the transfinite setting, we choose to consider a more involved recursion.

\begin{equation}\label{equation:boxBoxDefinition}
\begin{split}
\boxBox 0_T\phi \ \ \ &:= \Box_T \phi, \ \ \mbox{ and }\\
\boxBox{n+1}_T \phi  \ \ \        &:= \ \Box_T \phi \ \vee \ \exists \, \psi \  \bigvee_{0\leq m \leq n} \Big(\boxDiamond{m}_T \psi \ \wedge \ \Box (\boxDiamond{m}_T \psi \to \phi)\Big).
\end{split}
\end{equation}

We shall now prove that this provability notion $\boxBox n_T$ is actually \emph{provably} very similar to that of $\trueBox n_T$ for any natural number $n$. 

\begin{proposition}\label{theorem:boxBoxProvabilityProvablySimilarToTrueBoxProvability}
Let $T$ be a sound c.e.~theory extending \ea. We have for all 
$n \in \mathbb N$ that 
\begin{enumerate}
\item\label{item:boxBoxImpliesTrueBox:theorem:boxBoxProvabilityProvablySimilarToTrueBoxProvability}
$\ea \vdash \ \  \forall \varphi \ \big ( \boxBox n_T \varphi  \ \to \ \trueBox n_T \varphi \big)$;\\

\item\label{item:equivalentProvidedConsistency:theorem:boxBoxProvabilityProvablySimilarToTrueBoxProvability}
$\pa \vdash \ \  \trueDiamond{n}_T \top \ \ \to \ \  \forall \varphi \  \big ( \boxBox{n+1}_T \varphi \leftrightarrow  \trueBox{n+1}_T \varphi \big)$;\\

\item\label{item:equivalentUnderABox:theorem:boxBoxProvabilityProvablySimilarToTrueBoxProvability}
$\pa \vdash  \ \ \trueBox{n}_T \Big ( \ \  \forall \varphi \  \big ( \boxBox n_T \varphi \leftrightarrow  \trueBox n_T \varphi \big) \ \ \Big)$;\\

\item\label{item:onNequivalent:theorem:boxBoxProvabilityProvablySimilarToTrueBoxProvability}
$\mathbb N \models  \ \  \forall \varphi \  \big ( \boxBox n_T \varphi \leftrightarrow  \trueBox n_T \varphi \big) $.

\end{enumerate}
\end{proposition} 

\begin{proof}
For the sake of readability we shall omit the subscripts $T$ in this proof.

{\bf Item \ref{item:boxBoxImpliesTrueBox:theorem:boxBoxProvabilityProvablySimilarToTrueBoxProvability} :} 
This direction is easy, since by induction on $n$ we see that each $\boxDiamond n_T \psi$ is of complexity $\Pi_{n+1}$.


{\bf Item \ref{item:equivalentProvidedConsistency:theorem:boxBoxProvabilityProvablySimilarToTrueBoxProvability}:} We reason in \pa, assume $\trueDiamond{n} \top$ and pick $\varphi$ arbitrary. By Item \ref{item:boxBoxImpliesTrueBox:theorem:boxBoxProvabilityProvablySimilarToTrueBoxProvability} we only need to prove that $\trueBox{n+1}\varphi \to \boxBox{n+1}\varphi$. 

Thus, we suppose that $\trueBox{n+1}\varphi$ so that for some $\pi$ we have ${\sf True}_{\Pi_{n+1}}(\pi)$ and $\Box \big({\sf True}_{\Pi_{n+1}}(\pi) \to \varphi \big)$. Using the (formalized) FGH theorem we now pick $\rho$ so that ${\sf True}_{\Pi_{n+1}}(\pi) \wedge \trueDiamond n \top \leftrightarrow \trueDiamond n \rho$. Since we work under the assumption of $\trueDiamond n \top$ we thus have $\trueDiamond n \rho$. 

Clearly, since we know that $\Box \big({\sf True}_{\Pi_{n+1}}(\pi) \to \varphi \big)$ we also have the weaker  $\Box \big(\trueDiamond n \rho \to \varphi \big)$. Thus we have $\trueDiamond n \rho \wedge \Box \big(\trueDiamond n \rho \to \varphi \big)$ so that by applying twice the induction hypothesis we get
$\boxDiamond n \rho \wedge \Box \big(\boxDiamond n \rho \to \varphi \big)$ and by definition $\boxBox {n+1} \varphi$ as was to be shown.


{\bf Item \ref{item:equivalentUnderABox:theorem:boxBoxProvabilityProvablySimilarToTrueBoxProvability}:} For $n=0$ the statement holds by definition and for $n+1$, the statement follows from the previous item since $\ea \vdash \trueBox {n+1} \trueDiamond n \top$.

{\bf Item \ref{item:onNequivalent:theorem:boxBoxProvabilityProvablySimilarToTrueBoxProvability}:} follows directly from the previous from the soundness of $\ea$ and $T$.
\end{proof}

We shall see below that $\boxBox n_T$ provability is a very decent provability notion.

\begin{theorem}\label{theorem:recursionIsFineForFiniteNumbers}
Let $T$ be a c.e.~theory. We have for all 
$A \subseteq \mathbb N$ that the following are equivalent
\begin{enumerate}

\item\label{item:TuringReducible:theorem:recursionIsFineForFiniteNumbers}
$A$ is c.e.~in $\emptyset^{(n)}$;

\item\label{item:OneToOneReducible:theorem:recursionIsFineForFiniteNumbers}
$A$ is many-one reducible to $\emptyset^{(n+1)}$;

\item\label{item:BoxDefinable:theorem:recursionIsFineForFiniteNumbers}
$A$ is definable on the standard model by a formula of the form $\boxBox{n}_T \rho(\dot x)$.
\end{enumerate}
\end{theorem}

\begin{proof}
This is a direct consequence of Lemma \ref{theorem:TrueProvabilityTiesUpWithArithmeticalHierarchy} and a minor generalization of Proposition \ref{theorem:boxBoxProvabilityProvablySimilarToTrueBoxProvability}.
\end{proof}

As a consequence of this theorem we see that $\boxBox n_T$ sentences are closed under disjunctions and conjunctions.

\begin{corollary}\label{theorem:finiteboxBoxFormulasClosedUnderDisjunctions}
Let $T$ be a c.e.~theory extending \ea and let $n\in \mathbb N$. For each formulas $\varphi, \psi$ there is some $\sigma$ so that 
\[
T\vdash (\boxBox{n}_T \varphi \ \vee \boxBox{n}_T \psi) \ \leftrightarrow \ \boxBox{n}_T \sigma .
\]
\end{corollary}

\begin{proof}
Immediate from Theorem \ref{theorem:recursionIsFineForFiniteNumbers} since c.e.~sets (with or without oracles) are closed under both conjunctions and disjunctions. 
\end{proof}

As a matter of fact, it turns out that this corollary can be formalized being one of the corner stones in a proof to the extent that the provability logic concerning the $[n]^{\Box}$ predicates is nice. In particular, let \gl be the normal modal logic axiomatized by $\Box (A \to B) \to (\Box A \to \Box B)$ and $\Box (\Box A \to A) \to \Box A$ and all propositional tautologies with rules Modus Ponens and Necessitation: $\frac{A}{\Box A}$. Then $\glp_\omega$ is the polymodal logic axiomatized by \gl for each modality $[n]$ together with the schemas $[n]A \to [n+1]A$ and $\la n\ra A \to [n+1]\la n\ra A$. 

\begin{theorem}\label{theorem:glpOmegaSoundUnderBoxBoxInterpretation}
Let $T$ be any c.e.~theory extending \ea. The logic $\glp_\omega$ is sound w.r.t.~$\pa$ if we interpret each $[n]$-modality as $\boxBox{n}_T$.
\end{theorem}

The proof proceeds by a straightforward induction on $n$ considering the logics $\glp_{n+1}$ that have only modalities up to $[n]$.

\section*{Acknowledgements}
I would like to thank Lev Beklemishev, Ramon Jansana, Stephen Simpson and Albert Visser for encouragement and fruitful discussions/suggestions. 
%


\bibliographystyle{plain}
\bibliography{References}

\end{document}

%% file: Preamble.tex
\usepackage{amsmath,amssymb,amsthm,units,stmaryrd}
\usepackage{qtree,bussproofs}

\usepackage{yfonts}
\usepackage{units,stmaryrd}
\usepackage{color}
\usepackage{graphicx}
\usepackage[all]{xy}
\newtheorem{theorem}{Theorem}[section]

\newtheorem{definition}[theorem]{Definition}
\newtheorem{lemma}[theorem]{Lemma}
\newtheorem{corollary}[theorem]{Corollary}
\newtheorem{proposition}[theorem]{Proposition}

\newtheorem{fact}[theorem]{Fact}

\newcommand{\logic}[1]{{\ensuremath{\mathbf{#1}}}\xspace}

\newcommand{\glp}{{\ensuremath{\mathsf{GLP}}}\xspace}

\usepackage{xspace}

\newcommand{\pa}{\ensuremath{{\mathrm{PA}}}\xspace}

\newcommand{\gl}{\logic{GL}}

\newcommand{\ea}{\ensuremath{{\mathrm{EA}}}\xspace}

\newcommand{\la}{\langle}
\newcommand{\ra}{\rangle}

\def\fmodels{\xymatrix{
\ar@{|=}[r]^{<\omega}&
}
}
\def\nmodels{\xymatrix{
\ar@{|=}[r]^{N}&
}
}
\def\<{\left <}

\def\>{\right >}

\DeclareSymbolFont{AMSb}{U}{msb}{m}{n}
\DeclareMathSymbol{\N}{\mathbin}{AMSb}{"4E}
\DeclareMathSymbol{\Z}{\mathbin}{AMSb}{"5A}
\DeclareMathSymbol{\R}{\mathbin}{AMSb}{"52}
\DeclareMathSymbol{\Q}{\mathbin}{AMSb}{"51}
\DeclareMathSymbol{\I}{\mathbin}{AMSb}{"49}
\DeclareMathSymbol{\C}{\mathbin}{AMSb}{"43}

\newcommand{\bt}{\begin{theorem}}
\newcommand{\et}{\end{theorem}}
\newcommand{\bl}{\begin{lemma}}
\newcommand{\el}{\end{lemma}}

\newcommand\remove[1]{}

